\documentstyle[11pt]{article}
\begin{document}

\title{{\bf Gravity and the Noncommutative Residue for Manifolds with
Boundary}
\thanks{Supported by Science Foundation for Young Teachers of Northeast Normal University (No. 20060102)
 }}
\author{ Yong Wang\\
\thanks{also at School of Mathematics and Statistics, Northeast Normal University}
{\scriptsize \it  Center of Mathematical Sciences Zhejiang
University Hangzhou Zhejiang 310027, China ;}\\
{\scriptsize \it E-mail: wangy581@nenu.edu.cn}}

\date{}
\maketitle

\noindent {\bf Abstract}~~ We prove a Kastler-Kalau-Walze type
theorem for the Dirac operator and the signature operator for
$3,4$-dimensional manifolds with boundary. As a corollary, we give
two kinds of operator theoretic explanations of the gravitational
action in the case of $4$-dimensional manifolds with flat boundary.
\\
\noindent{\bf Subj. Class.:}\quad Noncommutative global analysis;
Noncommutative differential geometry.
\section{Introduction}
\quad The noncommutative residue found in [Gu] and [Wo] plays a
prominent role in noncommutative geometry. In [C1],Connes used the
noncommutative residue to derive a conformal 4-dimensional Polyakov
action analogy. In [C2], Connes proved that the noncommutative
residue on a compact manifold $M$ coincided with the Dixmier's trace
on pseudodifferential operators of order $-{\rm {dim}}M$. Several
years ago, Connes made a challenging observation that the
noncommutative residue of the square of the inverse of the Dirac
operator was proportional to the Einstein-Hilbert action, which we
call the Kastler-Kalau-Walze theorem. In [K], Kastler gave a
brute-force proof of this theorem. In [KW], Kalau and Walze proved
this theorem in the normal coordinates system simultaneously. In
[A], Ackermann gave a note on a new proof of this theorem
by means of the heat kernel expansion.\\
\indent On the other hand, Fedosov et al defined a noncommutative
residue on Boutet de Monvel's algebra and proved that it was a
unique continuous trace in [FGLS]. In [S], Schrohe gave the relation
between the Dixmier trace and the noncommutative residue for
manifolds with boundary. In [Wa1] and [Wa2], we generalized some
results in [C1] and [U] to the case of manifolds with boundary . In
[H], the gravitational action for manifolds with boundary was worked
out (also see [B]). The motivation of this paper is to give an
operator theoretic explanation of the gravitational action for
manifolds with boundary and prove a Kastler-Kalau-Walze type theorem
for manifolds with
boundary.\\
\indent Let us recall  the Kastler-Kalau-Walze theorem in [K],
[KW],[A]. Let $M$  be a $4$-dimensional oriented  spin manifold (it
holds
 for any even dimensional manifolds) and $D$ be the associated Dirac operator on the spinor bundle
$S(TM)$. Let $s$ be the scalar curvature and Wres denote  the
noncommutative residue (see [Wo],[FGV]). Then the
Kastler-Kalau-Walze theorem gives a spectral explanation of the
gravitational action, it says that there exists a constant $c_0$,
such that
$${\rm  Wres}(D^{-2})=c_0\int_Ms{\rm dvol}_M.\eqno(1.1)$$
For an oriented spin manifold  $M^4$ with boundary $\partial M$, we
use $ \widetilde{{\rm Wres}}[(\pi^+\widehat{D}^{-1})^2]$ instead of
${\rm Wres}(D^{-2})$ where $\widehat{D}$ is the Dirac operator on an
open neighborhood $\widehat{M}$ of $M$ and we still write $D$
instead of $\widehat{D}$ in this paper (for definition of
$\widehat{D}$ see Section 2) . Here $ \widetilde{{\rm Wres}}$
denotes the noncommutative residue for manifolds with boundary of
[FGLS] and $\pi^+D^{-1}$ is an element in Boutet de Monvel's algebra
(see [Wa1], Section 3).  By definitions of Boutet de Monvel's
algebra (see[S,p.11] or [Wa1,pp.5-6]), it is significant to consider
$$ \widetilde{{\rm Wres}}[(\pi^+\widehat{D}^{-1})^2]=\widetilde{{\rm Wres}}\left(\left[\begin{array}{lcr}
  \ \pi^+\widehat{D}^{-1} & 0 \\
    \  0 & 0
\end{array}\right]^2\right),\eqno(1.2)$$ \noindent which doesn't depend on the extension
  $\widehat{M}.$ By the composition formula
in Boutet de Monvel's algebra and the definition of $
\widetilde{{\rm Wres}}$ (see (2.4) and (2.6) in [Wa1]), $
\widetilde{{\rm Wres}}[(\pi^+D^{-1})^2]$ is the sum of two terms one
corresponding to interior and the other corresponding to boundary of
$M$. It is well known that (see [H]) that the gravitational action
for manifolds with boundary is also the sum of two terms from
interior and boundary of $M$.  So it is natural to hope to get the
gravitational action for manifolds with boundary by computing $
\widetilde{{\rm Wres}}[(\pi^+D^{-1})^2]$ . For simplicity, we assume
that the metric $g^M$ on $M$ has the following form near the
boundary,
$$ g^M=\frac{1}{h(x_n)}g^{\partial M}+dx_n^2,\eqno(1.3)$$
where $g^{\partial M}$ is the metric on  ${\partial M}$. $h(x_n)\in
C^{\infty}([0,1))=\{\widetilde{h}|_{[0,1)}|\widetilde{h}\in
C^{\infty}((-\varepsilon,1))\}$ for some $\varepsilon>0$ and
satisfies $h(x_n)>0,~h(0)=1$ where $x_n$ denotes the normal
directional coordinate.
 Through
computations, we find that the term from boundary
 which we expect to get vanishes, so $ \widetilde{{\rm
Wres}}[(\pi^+D^{-1})^2]$ is also proportional to $ \int_Ms{\rm
dvol}_M$.
 Fortunately, if we assume that $\partial M$ is flat, then we can define
$\int_{\partial M }{\rm res}_{1,1}(D^{-1},D^{-1})$ and
$\int_{\partial M} {\rm res}_{2,1}(D^{-1},D^{-1})$ (see Section 4)
and get that
 the gravitational  action for $\partial M$ is proportional to $\int_{\partial
 M}
{\rm res}_{1,1}(D^{-1},D^{-1})$ and $\int_{\partial M }{\rm
res}_{2,1}(D^{-1},D^{-1})$, which gives two kinds of operator
theoretic explanations
of the gravitational action for boundary. For general even dimensional manifolds with boundary,
we have no similar explanations for the gravitational action for boundary, even though for the flat
boundary (see Section 4).\\
\indent For odd dimensional manifolds without boundary, $ {\rm
Wres}(D^{-2})=0,$ so Kastler-Kalau-Walze Theorem isn't correct. But
for odd dimensional manifolds with boundary, in general $
\widetilde{{\rm Wres}}[(\pi^+\widehat{D}^{-1})^2]$ doesn't vanish
(similar to Section 5-7 in [Wa1]). In this paper we compute $
\widetilde{{\rm Wres}}[(\pi^+\widehat{D}^{-1})^2]$ explicitly for
$3$-dimensional
spin manifolds with boundary.\\
\indent This paper is organized as follows: In Section $2$, for
$4$-dimensional spin manifolds with boundary and the associated
Dirac operator $D$, we compute
 $\widetilde{{\rm
Wres}}[(\pi^+D^{-1})^2]$. In Section $3$, we compute
 $\widetilde{{\rm
Wres}}[(\pi^+D^{-1})^2]$ for $4$-dimensional oriented Riemannian
manifolds with boundary and the associated signature operator. Two
kinds of operator theoretic explanations of the gravitational action
for boundary in the case of $4$-dimensional manifolds with boundary
will be given in Section $4$.  In Section $5$, We compute $
\widetilde{{\rm Wres}}[(\pi^+\widehat{D}^{-1})^2]$ for
$3$-dimensional spin manifolds with boundary. In Appendix, the proof
of two facts in
Section $2$ will be given. \\

\section{The Dirac operator case}

\quad In this section, we compute
 $\widetilde{{\rm
Wres}}[(\pi^+D^{-1})^2]$ by the brute force way in [K] and the normal coordinates way in [KW].\\
\indent Let $M$ be a $n$-dimensional compact oriented spin manifold
with boundary $\partial M$ and the metric $g^M$ in (1.2). Let $n=4$,
but our some computations is correct for the general $n$. Let
$U\subset M$ be a collar neighborhood of $\partial M$ which is
diffeomorphic to $\partial M\times [0,1)$. By the definition of
$C^{\infty}([0,1))$ and $h>0$, there exists $\widetilde{h}\in
C^{\infty}((-\varepsilon,1))$ such that $\widetilde{h}|_{[0,1)}=h$
and $\widetilde{h}>0$ for some sufficiently small $\varepsilon>0$.
Then there exists a metric $\widehat{g}$ on
$\widehat{M}=M\cup_{\partial M}\partial M\times (-\varepsilon,0]$
which has the form on $U\cup_{\partial M}\partial M\times
(-\varepsilon,0]$
$$ \widehat{g}=\frac{1}{\widetilde{h}(x_n)}g^{\partial M}+dx_n^2,\eqno(2.1)$$
such that $\widehat{g}|_M=g.$ We fix a metric $\widehat{g}$ on the
$\widehat{M}$ such that $\widehat{g}|_M=g$. We can get the spin
structure on $\widehat{M}$ by extending the spin structure on $M.$
Let $D$ be the Dirac operator associated to $\widehat{g}$ on the
spinors bundle $S(T\widehat{M})$. We want to compute
 $\widetilde{{\rm Wres}}[(\pi^+D^{-1})^2]$ (for the related definitions, see [Wa1], Section 2, 3).
Let ${\bf S}~({\bf S}')$ be the unit sphere about $\xi~(\xi')$ and
$\sigma(\xi)~(\sigma(\xi'))$ be the corresponding canonical
$n-1~(n-2)$ volume form. Denote by $\sigma_l(A)$ the $l$- order
symbol of an operator $A$. By (2.4) and (2.6) in [Wa1], we get
$$\widetilde{{\rm Wres}}[(\pi^+D^{-1})^2]=\widetilde{{\rm
Wres}}[(\pi^+D^{-2})+L(D^{-1},D^{-1})]$$
$$=\int_M\int_{|\xi|=1}{\rm
trace}_{S(TM)}[\sigma_{-4}(D^{-2})]\sigma(\xi)dx +2\pi\int_
{\partial
M}\int_{|\xi'|=1}{\rm{tr}}_{S(TM)}[{\rm{tr}}(b_{-4})(x',\xi')]\sigma(\xi')dx',\eqno
(2.2)$$ where $b_{-4}$ is the $(-4)$-order symbol of
$L(D^{-1},D^{-1})$ which is called leftover term. By the formula
$(3.14)$ and $(3.15)$ in [Wa1] and $\pi'_{\xi_n}$ adding degree $1$
of the symbol and the $++$ parts vanishing after integration  with
respect to $\xi_n$ (see [FGLS] p. 23), then we have
$${\rm {tr}}(b_{-4})=
\sum_{j,k=0}^{\infty}\frac{(-i)^{j+k+1}}{(j+k+1)!}\pi'_{\xi_n}[\partial^j_{x_n}\partial^k_{\xi_n}
a^+(x',0,\xi',\xi_n)\circ'
\partial^{j+1}_{\xi_n}\partial^k_{x_n}a(x',0,\xi',\xi_n)]_{-4},\eqno(2.3)$$
\noindent where $a=\sigma(D^{-1})$ and $\pi'_{\xi_n}$ and
$a^+=\pi^+_{\xi_n}a$ defined by (2.1) and (2.2) in [Wa1]. By the
formula of p.740 line 2 in [Wa1], we get
$$\widetilde{{\rm
Wres}}[(\pi^+D^{-1})^2]=\int_M\int_{|\xi|=1}{\rm
trace}_{S(TM)}[\sigma_{-4}(D^{-2})]\sigma(\xi)dx+\int_{\partial
M}\Phi,\eqno(2.4)$$ \indent where
$$\Phi=\int_{|\xi'|=1}\int^{+\infty}_{-\infty}\sum^{\infty}_{j, k=0}
\sum\frac{(-i)^{|\alpha|+j+k+1}}{\alpha!(j+k+1)!}$$
$$\times {\rm trace}_{S(TM)}
[\partial^j_{x_n}\partial^\alpha_{\xi'}\partial^k_{\xi_n}
\sigma^+_{r}(D^{-1})(x',0,\xi',\xi_n)\times
\partial^\alpha_{x'}\partial^{j+1}_{\xi_n}\partial^k_{x_n}\sigma_{l}
(D^{-1})(x',0,\xi',\xi_n)]d\xi_n\sigma(\xi')dx',\eqno(2.5)$$
\noindent where the sum is taken over $
r-k-|\alpha|+l-j-1=-4,~~r,l\leq-1$. Since $[\sigma_{-4}(D^{-2})]|_M$
has the same expression as $\sigma_{-4}(D^{-2})$ in the case of
manifolds without boundary, so locally we can use the computations
in [K], [KW], [A], then we have
$$\int_M\int_{|\xi|=1}{\rm tr}[\sigma_{-4}(D^{-2})]\sigma(\xi)dx=-\frac{\Omega_4}{3}\int_Ms{\rm dvol}_M.\eqno(2.6)$$
where $\Omega_n=\frac{2\pi^{\frac{n}{2}}}{\Gamma(\frac{n}{2})}$. So we only need to compute $\int_{\partial M}\Phi$.\\
\indent Firstly, we compute the symbol $\sigma(D^{-1})$ of
$D^{-1}$. Recall the definition of the Dirac operator $D$ (see
[BGV], [Y]). Let $\nabla^L$ denote the Levi-civita connection
about $g^M$.
 In the local coordinates $\{x_i; 1\leq i\leq n\}$ and the fixed orthonormal frame $\{\widetilde{e_1},\cdots,\widetilde{e_n}\}$, the connection matrix $(\omega_{s,t})$
is defined by
$$\nabla^L(\widetilde{e_1},\cdots,\widetilde{e_n})= (\widetilde{e_1},\cdots,\widetilde{e_n})(\omega_{s,t}).\eqno(2.7)$$
$c(\widetilde{e_i})$ denotes the Clifford action. The Dirac
operator
$$D=\sum^n_{i=1}c(\widetilde{e_i})[\widetilde{e_i}
-\frac{1}{4}\sum_{s,t}\omega_{s,t}(\widetilde{e_i})c(\widetilde{e_s})c(\widetilde{e_t})].\eqno(2.8)$$
So we get,
$$\sigma_1(D)=\sqrt{-1}c(\xi); \sigma_0(D)
=-\frac{1}{4}\sum_{i,s,t}\omega_{s,t}(\widetilde{e_i})c(\widetilde{e_i})c(\widetilde{e_s})c(\widetilde{e_t}),\eqno(2.9)$$
where $\xi=\sum^n_{i=1}\xi_idx_i$ denotes the cotangent vector.
Write
$$D_x^{\alpha}=(-\sqrt{-1})^{|\alpha|}\partial_x^{\alpha};~\sigma(D)=p_1+p_0;
~\sigma(D^{-1})=\sum^{\infty}_{j=1}q_{-j}.\eqno(2.10)$$ By the
composition formula of psudodifferential operators, then we have
\begin{eqnarray*}
1=\sigma(D\circ D^{-1})&=&\sum_{\alpha}\frac{1}{\alpha!}\partial^{\alpha}_{\xi}[\sigma(D)]D^{\alpha}_{x}[\sigma(D^{-1})]\\
&=&(p_1+p_0)(q_{-1}+q_{-2}+q_{-3}+\cdots)\\
& &~~~+\sum_j(\partial_{\xi_j}p_1+\partial_{\xi_j}p_0)(
D_{x_j}q_{-1}+D_{x_j}q_{-2}+D_{x_j}q_{-3}+\cdots)\\
&=&p_1q_{-1}+(p_1q_{-2}+p_0q_{-1}+\sum_j\partial_{\xi_j}p_1D_{x_j}q_{-1})+\cdots,
\end{eqnarray*}
Thus, we get:
$$q_{-1}=p_1^{-1};~q_{-2}=-p_1^{-1}[p_0p_1^{-1}+\sum_j\partial_{\xi_j}p_1D_{x_j}(p_1^{-1})].\eqno(2.11)$$
By (2.9), (2.11) and direct computations, we have \\

\noindent{\bf Lemma 2.1}$$
q_{-1}=\frac{\sqrt{-1}c(\xi)}{|\xi|^2};~~~q_{-2}=\frac{c(\xi)p_0c(\xi)}{|\xi|^4}+\frac{c(\xi)}{|\xi|^6}\sum_jc(dx_j)
[\partial_{x_j}[c(\xi)]|\xi|^2-c(\xi)\partial_{x_j}(|\xi|^2)]\eqno(2.12)$$
\indent Since $\Phi$ is a global form on $\partial M$, so for any
fixed point $x_0\in\partial M$, we can choose the normal coordinates
$U$ of $x_0$ in $\partial M$ (not in $M$) and compute $\Phi(x_0)$ in
the coordinates $\widetilde{U}=U\times [0,1)\subset M$ and the
metric $\frac{1}{h(x_n)}g^{\partial M}+dx_n^2.$ The dual metric of
$g^M$ on $\widetilde{U}$ is ${h(x_n)}g^{\partial M}+dx_n^2.$  Write
$g^M_{ij}=g^M(\frac{\partial}{\partial x_i},\frac{\partial}{\partial
x_j});~ g_M^{ij}=g^M(dx_i,dx_j)$, then
$$[g^M_{i,j}]= \left[\begin{array}{lcr}
  \frac{1}{h(x_n)}[g_{i,j}^{\partial M}]  & 0  \\
   0  &  1
\end{array}\right];~~~
[g_M^{i,j}]= \left[\begin{array}{lcr}
  h(x_n)[g^{i,j}_{\partial M}]  & 0  \\
   0  &  1
\end{array}\right],\eqno(2.13)$$
and
$$\partial_{x_s}g_{ij}^{\partial M}(x_0)=0, 1\leq i,j\leq
n-1; ~~~g_{ij}^M(x_0)=\delta_{ij}.\eqno(2.14)$$ \indent Let $n=4$
and $\{e_1,\cdots,e_{n-1}\}$ be an orthonormal frame field in $U$
about $g^{\partial M}$ which is parallel along geodesics and
$e_i(x_0)=\frac{\partial}{\partial x_i}(x_0)$, then
$\{\widetilde{e_1}=\sqrt{h(x_n)}e_1,\cdots,\widetilde{e_{n-1}}=\sqrt{h(x_n)}e_{n-1},
\widetilde{e_n}=dx_n\}$ is the orthonormal frame field in
$\widetilde U$ about $g^M$. Locally $S(TM)|_{\widetilde {U}}\cong
\widetilde {U}\times\wedge^* _{\bf C}(\frac{n}{2}).$ Let
$\{f_1,\cdots,f_4\}$ be the orthonormal basis of $\wedge^* _{\bf
C}(\frac{n}{2}).$ Take a spin frame field $\sigma:~\widetilde
{U}\rightarrow {\rm Spin}(M)$ such that $\pi\sigma=
\{\widetilde{e_1},\cdots,\widetilde{e_n}\}$ where $\pi :~{\rm
Spin}(M)\rightarrow O(M)$ is a double covering, then
$\{[(\sigma,f_i)],~1\leq i\leq 4\}$ is an orthonormal frame of
$S(TM)|_{\widetilde {U}}.$ In the following, since the global form
$\Phi$ is independent of the choice of the local frame, so we can
compute ${\rm tr}_{S(TM)}$ in the frame $\{[(\sigma,f_i)],~1\leq
i\leq 4\}.$ Let $\{E_1,\cdots,E_n\}$ be the canonical basis of ${\bf
R}^n$ and $c(E_i)\in {\rm cl}_{\bf C}(n)\cong {\rm Hom}(\wedge^*
_{\bf C}(\frac{n}{2}),\wedge^* _{\bf C}(\frac{n}{2}))$ be the
Clifford action. By [Y], then
$$c( \widetilde{e_i})=[(\sigma,c(E_i))];~ c( \widetilde{e_i})[(\sigma,f_i)]=[(\sigma,c(E_i)f_i)];~
\frac{\partial}{\partial x_i}=[(\sigma,\frac{\partial}{\partial
x_i})],\eqno(2.15)$$
then we have $\frac{\partial}{\partial x_i}c( \widetilde{e_i})=0$ in the above frame.\\

\noindent{\bf Lemma 2.2}
$\partial_{x_j}(|\xi|_{g^M}^2)(x_0)=0,~{\rm if
}~j<n;~=h'(0)|\xi'|_{g^{\partial M}}^2,~{\rm if
}~j=n.~~~~~~(2.16)$\\
\noindent ~~~$\partial_{x_j}[c(\xi)](x_0)=0,~{\rm if
}~j<n;~=\partial_{x_n}[c(\xi')](x_0),~{\rm if
}~j=n,~~~~~~~(2.17)$\\
\noindent where $\xi=\xi'+\xi_ndx_n.$\\
 \noindent{\it Proof.}~~By the equality
  $\partial_{x_j}(|\xi|_{g^M}^2)(x_0)=\partial_{x_j}(h(x_n)g^{l,m}_{\partial M}(x')\xi_l\xi_m+\xi_n^2)(x_0)$ and
(2.14), then (2.16) is correct. By Lemma A.1 in Appendix, (2.17) is correct.\hfill$\Box$\\
\indent In order to compute $p_0(x_0)$, we need to compute $\omega_{s,t}(\widetilde{e_i})(x_0).$\\

\noindent{\bf Lemma 2.3}~~{\it When}~$
i<n,~\omega_{n,i}(\widetilde{e_i})(x_0)=\frac{1}{2}h'(0);$ {\it
and} $\omega_{i,n}(\widetilde{e_i})(x_0)=-\frac{1}{2}h'(0),$
{\it In other cases,} $\omega_{s,t}(\widetilde{e_i})(x_0)=0$\\
 \noindent{\it Proof.}~~See Appendix.\hfill$\Box$\\

\noindent{\bf Lemma 2.4}~~$p_0(x_0)=c_0c(dx_n),$ where $c_0=-\frac{3}{4}h'(0).$\\
\noindent{\bf Proof.}~~This comes from  (2.9),  Lemma 2.3 and
 the relation $c(\widetilde{e_i})c(\widetilde{e_j})+c(\widetilde{e_j})c(\widetilde{e_i})=-2\delta_{i.j}.$\hfill$\Box$\\
\indent Now we can compute $\Phi$, since the sum is taken over $
-r-l+k+j+|\alpha|=-3,~~r,l\leq-1,$ then we have the following five cases:\\

\noindent  {\bf case a)~I)}~$r=-1,~l=-1~k=j=0,~|\alpha|=1$\\

\noindent By (2.5), we get
$${\rm case~a)~I)}=-\int_{|\xi'|=1}\int^{+\infty}_{-\infty}\sum_{|\alpha|=1}
{\rm trace} [\partial^\alpha_{\xi'}\pi^+_{\xi_n}q_{-1}\times
\partial^\alpha_{x'}\partial_{\xi_n}q_{-1}](x_0)d\xi_n\sigma(\xi')dx',\eqno(2.17)$$
By Lemma 2.2, for $i<n$, then
$$\partial_{x_i}q_{-1}(x_0)=\partial_{x_i}\left(\frac{\sqrt{-1}c(\xi)}{|\xi|^2}\right)(x_0)=
\frac{\sqrt{-1}\partial_{x_i}[c(\xi)](x_0)}{|\xi|^2}
-\frac{\sqrt{-1}c(\xi)\partial_{x_i}(|\xi|^2)(x_0)}{|\xi|^4}=0,$$
\noindent so case a) I) vanishes.\\

\noindent  {\bf case a)~II)}~$r=-1,~l=-1~k=|\alpha|=0,~j=1$\\

\noindent By (2.5), we get
$${\rm case
a)~II)}=-\frac{1}{2}\int_{|\xi'|=1}\int^{+\infty}_{-\infty} {\rm
trace} [\partial_{x_n}\pi^+_{\xi_n}q_{-1}\times
\partial_{\xi_n}^2q_{-1}](x_0)d\xi_n\sigma(\xi')dx',\eqno(2.18)$$
\noindent By Lemma 2.1 and Lemma 2.2, we have\\
$$\partial^2_{\xi_n}q_{-1}=\sqrt{-1}\left(-\frac{6\xi_nc(dx_n)+2c(\xi')}
{|\xi|^4}+\frac{8\xi_n^2c(\xi)}{|\xi|^6}\right);\eqno(2.19)$$
$$\partial_{x_n}q_{-1}(x_0)
=\frac{\sqrt{-1}\partial_{x_n}c(\xi')(x_0)}{|\xi|^2}-\frac{\sqrt{-1}c(\xi)|\xi'|^2h'(0)}{|\xi|^4}.\eqno(2.20)$$
By (2.1) in [Wa1] and the Cauchy integral formula, then\\
\begin{eqnarray*}
\pi^+_{\xi_n}\left[\frac{c(\xi)}{|\xi|^4}\right](x_0)|_{|\xi'|=1}&=&\pi^+_{\xi_n}\left[\frac{c(\xi')+\xi_nc(dx_n)}{(1+\xi_n^2)^2}\right]\\
&=&\frac{1}{2\pi i}{\rm lim}_{u\rightarrow
0^-}\int_{\Gamma^+}\frac{\frac{c(\xi')+\eta_nc(dx_n)}{(\eta_n+i)^2(\xi_n+iu-\eta_n)}}
{(\eta_n-i)^2}d\eta_n\\
&=&\left[\frac{c(\xi')+\eta_nc(dx_n)}{(\eta_n+i)^2(\xi_n-\eta_n)}\right]^{(1)}|_{\eta_n=i}\\
&=&-\frac{ic(\xi')}{4(\xi_n-i)}-\frac{c(\xi')+ic(dx_n)}{4(\xi_n-i)^2}~~~~~~~~~~~~~~~~~~~(2.21)
\end{eqnarray*}
Similarly,
$$\pi^+_{\xi_n}\left[\frac{\sqrt{-1}\partial_{x_n}c(\xi')}{|\xi|^2}\right](x_0)|_{|\xi'|=1}
=\frac{\partial_{x_n}[c(\xi')](x_0)}{2(\xi_n-i)}.\eqno(2.22)$$
By (2.20), (2.21), (2.22), then\\
$$\pi^+_{\xi_n}\partial_{x_n}q_{-1}(x_0)|_{|\xi'|=1}=\frac{\partial_{x_n}[c(\xi')](x_0)}{2(\xi_n-i)}+\sqrt{-1}h'(0)
\left[\frac{ic(\xi')}{4(\xi_n-i)}+\frac{c(\xi')+ic(dx_n)}{4(\xi_n-i)^2}\right].\eqno(2.23)$$
\noindent By the relation of the Clifford action and ${\rm tr}{AB}={\rm tr }{BA}$, then we have the equalities:\\
$${\rm tr}[c(\xi')c(dx_n)]=0;~~{\rm tr}[c(dx_n)^2]=-4;~~{\rm tr}[c(\xi')^2](x_0)|_{|\xi'|=1}=-4;~~$$
$${\rm tr}[\partial_{x_n}c(\xi')c(dx_n)]=0;~~{\rm tr}[\partial_{x_n}c(\xi')c(\xi')](x_0)|_{|\xi'|=1}=-2h'(0).\eqno(2.24)$$
By (2.24) and direct computations ,we have
$$h'(0){\rm tr}\left\{\left[\frac{ic(\xi')}{4(\xi_n-i)}+\frac{c(\xi')+ic(dx_n)}{4(\xi_n-i)^2}\right]\times
\left[\frac{6\xi_nc(dx_n)+2c(\xi')}{(1+\xi_n^2)^2}-\frac{8\xi_n^2[c(\xi')+\xi_nc(dx_n)]}{(1+\xi_n^2)^3}\right]
\right\}(x_0)|_{|\xi'|=1}$$
$$=-4h'(0)\frac{-2i\xi_n^2-\xi_n+i}{(\xi_n-i)^4(\xi_n+i)^3}.\eqno(2.25)$$
Similarly, we have
$$-\sqrt{-1}{\rm
tr}\left\{\left[\frac{\partial_{x_n}[c(\xi')](x_0)}{2(\xi_n-i)}\right]
\times\left[\frac{6\xi_nc(dx_n)+2c(\xi')}{(1+\xi_n^2)^2}-\frac{8\xi_n^2[c(\xi')+\xi_nc(dx_n)]}
{(1+\xi_n^2)^3}\right]\right\}(x_0)|_{|\xi'|=1}$$
$$=-2\sqrt{-1}h'(0)\frac{3\xi_n^2-1}{(\xi_n-i)^4(\xi_n+i)^3}.\eqno(2.26)$$
By (2.19), (2.23), (2.25), (2.26), then\\
\begin{eqnarray*}
{\rm case~ a)~
II)}&=&-\int_{|\xi'|=1}\int^{+\infty}_{-\infty}\frac{ih'(0)(\xi_n-i)^2}
{(\xi_n-i)^4(\xi_n+i)^3}d\xi_n\sigma(\xi')dx'\\
&=&-ih'(0)\Omega_3\int_{\Gamma^+}\frac{1}{(\xi_n-i)^2(\xi_n+i)^3}d\xi_ndx'\\
&=&-ih'(0)\Omega_32\pi i[\frac{1}{(\xi_n+i)^3}]^{(1)}|_{\xi_n=i}dx'\\
&=&-\frac{3}{8}\pi h'(0)\Omega_3dx'.
\end{eqnarray*}

\noindent  {\bf case a)~III)}~$r=-1,~l=-1~j=|\alpha|=0,~k=1$\\

\noindent By (2.5), we get
$${\rm case~ a)~III)}=-\frac{1}{2}\int_{|\xi'|=1}\int^{+\infty}_{-\infty}
{\rm trace} [\partial_{\xi_n}\pi^+_{\xi_n}q_{-1}\times
\partial_{\xi_n}\partial_{x_n}q_{-1}](x_0)d\xi_n\sigma(\xi')dx',\eqno(2.27)$$
\noindent By Lemma 2.2, we have\\
$$\partial_{\xi_n}\partial_{x_n}q_{-1}(x_0)|_{|\xi'|=1}=-\sqrt{-1}h'(0)
\left[\frac{c(dx_n)}{|\xi|^4}-4\xi_n\frac{c(\xi')+\xi_nc(dx_n)}{|\xi|^6}\right]-
\frac{2\xi_n\sqrt{-1}\partial_{x_n}c(\xi')(x_0)}{|\xi|^4}.\eqno(2.28)$$
$$\partial_{\xi_n}\pi^+_{\xi_n}q_{-1}(x_0)|_{|\xi'|=1}=-\frac{c(\xi')+ic(dx_n)}{2(\xi_n-i)^2}.\eqno(2.29)$$
Similarly to (2.25), (2.26), we have\\
$${\rm tr}\left\{\frac{c(\xi')+ic(dx_n)}{2(\xi_n-i)^2}\times
\sqrt{-1}h'(0)\left[\frac{c(dx_n)}{|\xi|^4}-4\xi_n\frac{c(\xi')+\xi_nc(dx_n)}{|\xi|^6}\right]\right\}$$
$$=2h'(0)\frac{i-3\xi_n}{(\xi_n-i)^4(\xi_n+i)^3};\eqno(2.30)$$
\noindent and\\
$${\rm tr}\left[\frac{c(\xi')+ic(dx_n)}{2(\xi_n-i)^2}\times
\frac{2\xi_n\sqrt{-1}\partial_{x_n}c(\xi')(x_0)}{|\xi|^4}\right]
=-2h'(0)\sqrt{-1}\frac{\xi_n}{(\xi_n-i)^4(\xi_n+i)^2}.\eqno(2.31)$$
\noindent So we get case a) III)$=\frac{3}{8}\pi h'(0)\Omega_3dx'.$\\

\noindent  {\bf case b)}~$r=-2,~l=-1,~k=j=|\alpha|=0$\\

\noindent By (2.5), we get
$${\rm case~ b)}=-i\int_{|\xi'|=1}\int^{+\infty}_{-\infty}
{\rm trace} [\pi^+_{\xi_n}q_{-2}\times
\partial_{\xi_n}q_{-1}](x_0)d\xi_n\sigma(\xi')dx',\eqno(2.32)$$
\noindent By Lemma 2.1 and Lemma 2.2, we have\\
$$q_{-2}(x_0)=\frac{c(\xi)p_0(x_0)c(\xi)}{|\xi|^4}+\frac{c(\xi)}{|\xi|^6}c(dx_n)
[\partial_{x_n}[c(\xi')](x_0)|\xi|^2-c(\xi)h'(0)|\xi|^2_{\partial
M}].\eqno(2.33)$$ \noindent Then
$$\pi^+_{\xi_n}q_{-2}(x_0)|_{|\xi'|=1}
=\pi^+_{\xi_n}\left[\frac{c(\xi)p_0(x_0)c(\xi)+c(\xi)c(dx_n)\partial_{x_n}[c(\xi')](x_0)}{(1+\xi_n^2)^2}\right]$$
$$-h'(0)\pi^+_{\xi_n}\left[\frac{c(\xi)c(dx_n)c(\xi)}{(1+\xi_n)^3}\right]
:=B_1-B_2.\eqno (2.34)$$
\noindent Similarly to (2.21), we have\\
$$B_1=-\frac{A_1}{4(\xi_n-i)}-\frac{A_2}{4(\xi_n-i)^2},\eqno(2.35)$$
\noindent where
$$A_1=ic(\xi')p_0c(\xi')+ic(dx_n)p_0c(dx_n)+ic(\xi')c(dx_n)\partial_{x_n}[c(\xi')];$$
$$A_2=[c(\xi')+ic(dx_n)]p_0[c(\xi')+ic(dx_n)]+c(\xi')c(dx_n)\partial_{x_n}c(\xi')-i\partial_{x_n}[c(\xi')].\eqno(2.36)$$
\begin{eqnarray*}
B_2&=&h'(0)\pi_{\xi_n}^+\left[\frac{-\xi_n^2c(dx_n)^2-2\xi_nc(\xi')+c(dx_n)}{(1+\xi_n^2)^3}\right]\\
&=&\frac{h'(0)}{2}\left[\frac{-\eta^2_nc(dx_n)-2\eta_nc(\xi')+c(dx_n)}{(\eta_n+i)^3(\xi_n-\eta_n)}\right]^{(2)}|_{\eta_n=i}\\
&=&\frac{h'(0)}{2}\left[\frac{c(dx_n)}{4i(\xi_n-i)}+\frac{c(dx_n)-ic(\xi')}{8(\xi_n-i)^2}
+\frac{3\xi_n-7i}{8(\xi_n-i)^3}[ic(\xi')-c(dx_n)]\right].(2.37)
\end{eqnarray*}
$$\partial_{\xi_n}q_{-1}(x_0)|_{|\xi'|=1}=\sqrt{-1}
\left[\frac{c(dx_n)}{1+\xi_n^2}-\frac{2\xi_nc(\xi')+2\xi_n^2c(dx_n)}{(1+\xi_n^2)^2}\right].\eqno(2.38)$$
\noindent By (2.37), (2.38), we have\\
$${\rm tr }[B_2\times\partial_{\xi_n}q_{-1}(x_0)]|_{|\xi'|=1}
=\frac{\sqrt{-1}}{2}h'(0){\rm trace}$$
$$\left\{\{\left[\frac{1}{4i(\xi_n-i)}+\frac{1}{8(\xi_n-i)^2}-\frac{3\xi_n-7i}{8(\xi_n-i)^3}\right]c(dx_n)
+\left[\frac{-1}{8(\xi_n-i)^2}+\frac{3\xi_n-7i}{8(\xi_n-i)^3}\right]ic(\xi')\}\right.$$
$$\times\left.\{\left[\frac{1}{1+\xi_n^2}-\frac{2\xi_n^2}{(1+\xi_n^2)^2}\right]c(dx_n)-\frac{2\xi_n}
{(1+\xi_n^2)^2}c(\xi')\}\right\}$$
$$=\frac{\sqrt{-1}}{2}h'(0)\frac{-i\xi_n^2-\xi_n+4i}{4(\xi_n-i)^3(\xi_n+i)^2}{\rm tr}[{\rm id}].
~~~~~~~~~~~~~~~~~~~~~~~~~~~~~~\eqno(2.39)$$ \noindent Note that
$$B_1=\frac{-1}{4(\xi_n-i)^2}[(2+i\xi_n)c(\xi')p_0c(\xi')+i\xi_nc(dx_n)p_0c(dx_n)$$
$$+
(2+i\xi_n)c(\xi')c(dx_n)\partial_{x_n}c(\xi')+ic(dx_n)p_0c(\xi')
+ic(\xi')p_0c(dx_n)-i\partial_{x_n}c(\xi')].\eqno(2.40)$$ \noindent
By (2.24), (2.38), (2.40), Lemma 2.4 and ${\rm tr}(AB)={\rm
tr}(BA)$,
 considering for $i<n$ \\
\noindent $\int_{|\xi'|=1}\{\xi_{i_1}\xi_{i_2}\cdots\xi_{i_{2d+1}}\}\sigma(\xi')=0$, then\\
$${\rm tr }[B_1\times\partial_{\xi_n}q_{-1}(x_0)]|_{|\xi'|=1}=
\frac{-2ic_0}{(1+\xi_n^2)^2}+h'(0)\frac{\xi_n^2-i\xi_n-2}{2(\xi_n-i)(1+\xi_n^2)^2}.\eqno(2.41)$$
By (2.34), (2.39) and  (2.41), we have\\
$${\rm case~ b)}=-\Omega_3\int_{\Gamma_+}\frac{2c_0(\xi_n-i)+ih'(0)}{(\xi_n-i)^3(\xi_n+i)^2}d\xi_ndx'=
\frac{9}{8}\pi h'(0)\Omega_3dx'.\eqno(2.42)$$

\noindent {\bf  case c)}~$r=-1,~l=-2,~k=j=|\alpha|=0$\\

\noindent By (2.5), we get
$${\rm case~ c)}=-i\int_{|\xi'|=1}\int^{+\infty}_{-\infty}
{\rm trace} [\pi^+_{\xi_n}q_{-1}\times
\partial_{\xi_n}q_{-2}](x_0)d\xi_n\sigma(\xi')dx'.\eqno(2.43)$$
\noindent By
$$\pi^+_{\xi_n}q_{-1}(x_0)|_{|\xi'|=1}=\frac{c(\xi')+ic(dx_n)}{2(\xi_n-i)};\eqno(2.44)$$
$$\partial_{\xi_n}q_{-2}(x_0)|_{|\xi'|=1}=
\frac{1}{(1+\xi_n^2)^3}[(2\xi_n-2\xi_n^3)c(dx_n)p_0c(dx_n)+(1-3\xi_n^2)c(dx_n)p_0c(\xi')$$
$$+ (1-3\xi_n^2)c(\xi')p_0c(dx_n)-4\xi_nc(\xi')p_0c(\xi')
+(3\xi_n^2-1)\partial_{x_n}c(\xi')-4\xi_nc(\xi')c(dx_n)\partial_{x_n}c(\xi')$$
$$+2h'(0)c(\xi')+2h'(0)\xi_nc(dx_n)]
+6\xi_nh'(0)\frac{c(\xi)c(dx_n)c(\xi)}{(1+\xi^2_n)^4},\eqno(2.45)$$
\noindent then similarly to computations of the case b), we have
$${\rm trace}
[\pi^+_{\xi_n}q_{-1}\times
\partial_{\xi_n}q_{-2}](x_0)|_{|\xi'|=1}=\frac{3h'(0)(i\xi_n^2
+\xi_n-2i)}{(\xi_n-i)^3(\xi_n+i)^3}+\frac{12h'(0)i\xi_n}{(\xi_n-i)^3(\xi_n+i)^4}.\eqno(2.46)$$
\noindent So ${\rm case~ c)}=-\frac{9}{8}\pi h'(0)\Omega_3dx'.$ Now $\Phi$ is the sum of the cases a), b) and c), so is zero. Then we get\\

\noindent {\bf Theorem 2.5}~~{\it  Let  $M$ be a $4$-dimensional
compact spin manifold with the boundary $\partial M$ and the metric
$g^M$ as above and $D$ be the Dirac operator on $\widehat{M}$, then}
$$ \widetilde{{\rm
Wres}}[(\pi^+D^{-1})^2]=-\frac{\Omega_4}{3}\int_Ms{\rm
dvol}_M.\eqno(2.47)$$

 \noindent {\bf Remark 2.6}~ Since (2.4) is correct for any
 dimensional manifolds with boundary, we conjecture that Theorem 2.5 is also
 correct for any even dimensional manifolds with boundary. But our
 computations way maybe isn't valid for general even dimensional manifolds with
 boundary. When the dimension becomes larger and larger, the terms
 which we need to compute becomes more and more. Maybe the way in
 [GS] is valid for any even dimensional manifolds with boundary.\\

\section{The signature operator case}

\quad Let $M$ be a $4$-dimensional compact oriented Riemannian
manifold with boundary $\partial M$ and the metric in Section 2.
$D=d+\delta:~\wedge^*(T^*M)\rightarrow \wedge^*(T^*M)$ is the
signature operator. Take the coordinates and the orthonormal frame
as in Section 2.
 Let $\epsilon (\widetilde{e_j*} ),~\iota (\widetilde{e_j*} )$ be the exterior and interior multiplications respectively. Write
$$c(\widetilde{e_j})=\epsilon (\widetilde{e_j*} )-\iota (\widetilde{e_j*} );~~
\bar{c}(\widetilde{e_j})=\epsilon (\widetilde{e_j*} )+\iota
(\widetilde{e_j*} ).\eqno(3.1)$$ \noindent We'll compute ${\rm
tr}_{\wedge^*(T^*M)}$ in the frame $\{dx_{i_1}\wedge\cdots\wedge
dx_{i_k}|~1\leq i_1<\cdots<i_k\leq 4\}.$ By [Y], we have
$$D=d+\delta=\sum^n_{i=1}c(\widetilde{e_i})[\widetilde{e_i}+\frac{1}{4}\sum_{s,t}\omega_{s,t}
(\widetilde{e_i})[\bar{c}(\widetilde{e_s})\bar{c}(\widetilde{e_t})-c(\widetilde{e_s})c(\widetilde{e_t})]].\eqno(3.2)$$
So
$$p_1=\sigma_1(d+\delta)=\sqrt{-1}c(\xi);~p_0=\sigma_0(d+\delta)=
\frac{1}{4}\sum_{i,s,t}\omega_{s,t}
(\widetilde{e_i})c(\widetilde{e_i})[\bar{c}(\widetilde{e_s})\bar{c}(\widetilde{e_t})
-c(\widetilde{e_s})c(\widetilde{e_t})].\eqno(3.3)$$ \noindent Lemmas
2.1-2.3 are also correct, by Lemma 2.3, then
$$p_0(x_0)=\widetilde{p_0}(x_0)-\frac{3}{4}h'(0)c(dx_n);~
\widetilde{p_0}(x_0)=\frac{1}{4}h'(0)\sum^{n-1}_{i=1}c(\widetilde{e_i})\bar{c}(\widetilde{e_n})\bar{c}(\widetilde{e_i})(x_0).\eqno(3.4)$$
\noindent For the signature operator case,
$$ {\rm tr}[{\rm id}]=16;~~{\rm tr}[c(\xi')\partial_{x_n}c(\xi')](x_0)|_{|\xi'|=1}=-8h'(0);\eqno(3.5)$$
$${\rm tr}[c(\xi')p_0c(\xi')c(dx_n)](x_0)
={\rm tr}[p_0c(\xi')c(dx_n)c(\xi')](x_0)=|\xi'|^2{\rm
tr}[p_0c(dx_n)].\eqno(3.6)$$
\begin{eqnarray*}
c(dx_n)\widetilde{p_0}(x_0)&=&-\frac{1}{4}h'(0)\sum^{n-1}_{i=1}c(e_i)\bar{c}(e_i)c(e_n)\bar{c}(e_n)\\
&=&-\frac{1}{4}h'(0)\sum^{n-1}_{i=1}[\epsilon ({e_i*} )\iota
({e_i*} )-\iota(e_i*)\epsilon(e_i*)][\epsilon ({e_n*} )\iota
({e_n*} )-\iota(e_n*)\epsilon(e_n*)]
\end{eqnarray*}
\noindent By Theorem 4.3 in [U], then
$${\rm tr}_{\wedge^m(T^*M)} \{[\epsilon ({e_i*} )\iota ({e_i*} )-\iota(e_i*)\epsilon(e_i*)]
[\epsilon ({e_n*} )\iota ({e_n*} )-\iota(e_n*)\epsilon(e_n*)]\}$$
$$=a_{n,m}<e_i*,e_n*>^2+b_{n,m}|e_i*|^2|e_n*|^2=b_{n,m},\eqno(3.7)$$
\noindent where $b_{4,m}=\left(\begin{array}{lcr}
  \ \ 2 \\
    \  m-2
\end{array}\right)+\left(\begin{array}{lcr}
  \ \ 2 \\
    \  m
\end{array}\right)-2\left(\begin{array}{lcr}
  \ \ 2 \\
    \  m-1
\end{array}\right).$  By (3.7), then
$${\rm tr}_{\wedge^*(T^*M)} \{[\epsilon ({e_i*} )\iota ({e_i*} )-\iota(e_i*)\epsilon(e_i*)][\epsilon ({e_n*} )\iota ({e_n*} )-\iota(e_n*)\epsilon(e_n*)]\}
=\sum_{m=0}^4b_{4,m}=0.$$ \noindent Then
$${\rm
tr}_{\wedge^*(T^*M)}[c(dx_n)\widetilde{p_0}(x_0)]=0.\eqno(3.8)$$
\noindent By  (3.4), (3.5), (3.6) and (3.8), then $\Phi_{\rm{sig}}=4\Phi_{\rm {Dirac}}=0.$ So we get\\

\noindent {\bf Theorem 3.1}~~{\it  Let  $M$ be a $4$-dimensional
compact oriented Riemaniann manifold with the boundary $\partial M$
and the metric $g^M$ as above and $D$ be the signature operator on
$\widehat{M}$, then}
$$ \widetilde{{\rm
Wres}}[(\pi^+D^{-1})^2]=\frac{8\Omega_4}{3}\int_Ms{\rm
dvol}_M.\eqno(3.9)$$

\section{The gravitational action for $4$-dimensional manifolds with boundary}

 \quad Firstly, we recall the Einstein-Hilbert action for manifolds with boundary (see [H] or [B]),
$$I_{\rm Gr}=\frac{1}{16\pi}\int_Ms{\rm dvol}_M+2\int_{\partial M}K{\rm dvol}_{\partial_M}:=I_{\rm {Gr,i}}+I_{\rm {Gr,b}},\eqno(4.1)$$
\noindent  where
$$ K=\sum_{1\leq i,j\leq {n-1}}K_{i,j}g_{\partial M}^{i,j};~~K_{i,j}=-\Gamma^n_{i,j},\eqno(4.2)$$
\noindent and $K_{i,j}$ is the second fundamental form, or extrinsic
curvature. Take the metric in Section 2, then by Lemma A.2,
$K_{i,j}(x_0)=-\Gamma^n_{i,j}(x_0)=-\frac{1}{2}h'(0),$ when $i=j<n$,
otherwise is zero. For $n=4$, then
$$K(x_0)=\sum_{i,j}K_{i.j}(x_0)g_{\partial M}^{i,j}(x_0)=\sum_{i=1}^3K_{i,i}(x_0)=-\frac{3}{2}h'(0),$$
\noindent So
 $$I_{\rm {Gr,b}}=-3h'(0){\rm Vol}_{\partial M}.\eqno(4.3)$$
\indent Let $M$ be $4$-dimensional manifolds with boundary and
$P,P'$ be two pseudodifferential operators with transmission
property (see [Wa1] or [RS]) on $\widehat M$. By (2.4) in [Wa1], we
have $$\pi^+P\circ\pi^+P'=\pi^+(PP')+L(P,P')\eqno(4.4)$$ and
$L(P,P')$ is leftover term which represents the difference between
the composition $\pi^+P\circ\pi^+P'$ in Boutet de Monvel algebra and
the composition $PP'$ in the classical pseudodifferential operators
algebra. By (2.5), we define locally
 $${\rm
res}_{1,1}(P,P'):=-\frac{1}{2}\int_{|\xi'|=1}\int^{+\infty}_{-\infty}
{\rm trace} [\partial_{x_n}\pi^+_{\xi_n}\sigma_{-1}(P)\times
\partial_{\xi_n}^2\sigma_{-1}(P')]d\xi_n\sigma(\xi')dx';\eqno(4.5)$$
$${\rm res}_{2,1}(P,P'):=-i\int_{|\xi'|=1}\int^{+\infty}_{-\infty}
{\rm trace} [\pi^+_{\xi_n}\sigma_{-2}(P)\times
\partial_{\xi_n}\sigma_{-1}(P')]d\xi_n\sigma(\xi')dx'.\eqno(4.6)$$
\noindent Thus they represent the difference between the composition
$\pi^+P\circ\pi^+P'$ in Boutet de Monvel algebra and the composition
$PP'$ in the classical pseudodifferential operators algebra
partially and
$$ {\rm case~ a)~ II)}={\rm res}_{1,1}(D^{-1},D^{-1});~{\rm case~ b)}={\rm res}_{2,1}(D^{-1},D^{-1}).\eqno(4.7)$$
\noindent Now, we assume $\partial M$ is flat , then
$\{dx_i=e_i\},~g^{\partial M}_{i,j}=\delta_{i,j},~\partial
_{x_s}g^{\partial M}_{i,j}=0$. So ${\rm res}_{1,1}(D^{-1},D^{-1})$
and ${\rm res}_{2,1}(D^{-1},D^{-1})$ are two global forms locally
defined by the aboved oriented orthonormal basis $\{dx_i\}$. By case
a) II) and case b),
then we have:\\

\noindent {\bf Theorem 4.1}~~{\it  Let  $M$ be a $4$-dimensional
compact spin manifold with the boundary $\partial M$ and the metric
$g^M$ as above and $D$ be the Dirac operator on $\widehat{M}$.
Assume $\partial M$ is flat, then}
$$\int_{\partial M}{\rm res}_{1,1}(D^{-1},D^{-1})=\frac{\pi}{8}\Omega_3I_{\rm {Gr,b}};\eqno(4.8)$$
$$\int_{\partial M}{\rm res}_{2,1}(D^{-1},D^{-1})=-\frac{3\pi}{8}\Omega_3I_{\rm {Gr,b}}.\eqno(4.9)$$\\

\noindent {\bf Theorem 4.2}~~{\it  Let  $M$ be a $4$-dimensional
compact oriented Riemaniann manifold with the boundary $\partial M$
and the metric $g^M$ as above and $D$ be the signature operator on
$\widehat{M}$. Assume $\partial M$ is flat, then}
$$\int_{\partial M}{\rm res}_{1,1}(D^{-1},D^{-1})=\frac{\pi}{2}\Omega_3I_{\rm {Gr,b}};\eqno(4.10)$$
$$\int_{\partial M}{\rm res}_{2,1}(D^{-1},D^{-1})=-\frac{3\pi}{2}\Omega_3I_{\rm {Gr,b}}.\eqno(4.11)$$\\

\noindent {\bf Remark 4.3}~ We take $N$ is a flat 3-dimensional
oriented Riemannian manifold and $M=N\times [0,1]$, then $\partial
M=N\oplus N.$ Let $g^M=\frac{1}{h(x_n)}g^N+dx_n^2$, where
$h(x_n)=1-x_n(x_n-1)>0$ for $x_n\in [0,1]$ and $h(0)=h(1)=1.$ The
$(M,g^M)$ satisfies conditions in Theorem 4.2. Similar construction
is correct for Theorem 4.1. When $\partial M$ is not connected, we
still define the
noncommutative residue with the loss of the unique property.\\

\noindent {\bf Remark 4.4}~ Considering (2.5). when the dimension
increases, the degree of the derivative of $h(x_n)$ in $\Phi$ will
increase. So the theorems 4.1 and 4.2 aren't correct for any even
dimensional
manifolds.\\

 \noindent {\bf Remark 4.5}~ The reason
that the term from boundary does not appear is perhaps that we
ignore boundary conditions. We hope to compute the noncommutative
residue $\widetilde{{\rm Wres}}[(\pi^+D^{-1})^2]$ under certain
boundary conditions to get the term from boundary in the future.
Grubb and Schrohe got the noncommutative residue for manifolds with
boundary through asymptotic expansions in [GS]. Another problem is
to compute
 $\widetilde{{\rm Wres}}[(\pi^+D^{-1})^2]$ by asymptotic
 expansions.\\

\section{Computations of  $ \widetilde{{\rm
Wres}}[(\pi^+\widehat{D}^{-1})^2]$ for $3$-dimensional spin
manifolds with boundary}

\quad For an odd dimensional manifolds with boundary, as in Section
5-7 in [Wa1], we have the formula
$$\widetilde{{\rm
Wres}}[(\pi^+D^{-1})^2]=\int_{\partial M}\Phi.\eqno(5.1)$$ \noindent
When $n=3$, then in (2.5), $ r-k-|\alpha|+l-j-1=-3,~~r,l\leq-1$, so
we get $r=l=-1,~k=|\alpha|=j=0,$ then
$$\Phi=\int_{|\xi'|=1}\int^{+\infty}_{-\infty}
 {\rm trace}_{S(TM)}
[ \sigma^+_{-1}(D^{-1})(x',0,\xi',\xi_n)\times
\partial_{\xi_n}\sigma_{-1}
(D^{-1})(x',0,\xi',\xi_n)]d\xi_3\sigma(\xi')dx'.\eqno(5.2)$$
\noindent By Lemma 2.1, then similar to (2.21), we have
$$\sigma^+_{-1}(D^{-1})|_{|\xi'|=1}=\frac{\sqrt{-1}[c(\xi')+ic(dx_n)]}{2i(\xi_n-i)};\eqno(5.3)$$
$$\partial_{\xi_n}\sigma_{-1}(D^{-1})|_{|\xi'|=1}=\frac{\sqrt{-1}c(dx_n)}{1+\xi_n^2}
-\frac{2\sqrt{-1}\xi_nc(\xi)}{(1+\xi_n^2)^2}.\eqno(5.4)$$ \noindent
For $n=3$, we take the coordinates as in Section 2.
 Locally $S(TM)|_{\widetilde {U}}\cong
\widetilde {U}\times\wedge^{{\rm even}} _{\bf C}(2).$ Let
$\{\widetilde{f_1},\widetilde{f_2}\}$ be an orthonormal basis of
$\wedge^{{\rm even}} _{\bf C}(2)$ and we will compute the trace
under this basis. Similarly to (2.24), we have
$${\rm tr}[c(\xi')c(dx_3)]=0;~~{\rm tr}[c(dx_3)^2]=-2;~~{\rm tr}[c(\xi')^2](x_0)|_{|\xi'|=1}=-2
\eqno(5.5)$$ Then by (5.3) (5.4) and (5.5), we get
$${\rm trace} [
\sigma^+_{-1}(D^{-1})\times
\partial_{\xi_n}\sigma_{-1}
(D^{-1})](x_0)|_{|\xi'|=1}\\
=-\frac{1}{(\xi_n+i)^2(\xi_n-i)}.\eqno(5.6)$$
 \noindent By (5.2) and (5.6)
and the Cauchy integral formula, we get
$$\Phi=\frac{i\pi}{2}\Omega_2{\rm vol}_{\partial M}=i\pi^2{\rm vol}_{\partial
M}.\eqno(5.7)$$ \noindent Here ${\rm vol}_{\partial M}$ denotes the
canonical volume form of ${\partial M}$.\\

\noindent {\bf Theorem 5.1}~~{\it Let $M$ be a $3$-dimensional
compact spin manifold with the boundary $\partial M$ and the metric
$g^M$ as in Section 2 and $D$ be the Dirac operator on $\widehat{M}$
, then}
$$ \widetilde{{\rm
Wres}}[(\pi^+D^{-1})^2]=i\pi^2{\rm Vol}_{\partial M},\eqno(5.8)$$
{\it where ${\rm Vol}_{\partial M}$ denotes the canonical volume of
${\partial M}.$}\\

 \noindent {\bf Remark 5.2}~By Theorem 5.1, we know
that $\widetilde{{\rm Wres}}[(\pi^+D^{-1})^2]$ isn't proportional to
the gravitational action for boundary for $3$-dimensional manifolds
with boundary. By the same reason as in Remark 4.4, we know that
$\widetilde{{\rm Wres}}[(\pi^+D^{-1})^2]$ isn't proportional to the
gravitational action for boundary for any odd dimensional manifolds
with boundary.\\

\noindent {\bf Appendix}\\

\quad In this appendix, we will prove some facts used in Lemma 2.2 and Lemma 2.3.\\

\noindent {\bf Lemma A.1} $$\partial_{x_l}c(dx_j)(x_0)=0, {\it
when}~l<n;~\partial_{x_l}c(dx_n)=0$$
\noindent{\bf  Proof.}  The
fundamental setup is as in Section 2. Write $<\partial
_{x_s},e_i>_{g^{\partial M}}=H_{i,s}$, then by [Y] or [BGV],
$\partial_{x_j}H_{i,s}(x_0)=0$. Define $dx_j^*\in TM|_{\widetilde
U}$ by $<dx_j^*,v>=(dx_j,v)$ for $v\in TM.$ For $j<n$,
\begin{eqnarray*}
c(dx_j)&=&c(dx_j^*)=c(\sum_i<dx_j^*,\widetilde{e_i}>\widetilde{e_i})\\
&=&\sum_{i,s}g^{s,j}<\partial_{x_s},\widetilde{e_i}>_{g^M}c(\widetilde
{e_i})=\sum_{1\leq
i,s<n}\frac{1}{\sqrt{h(x_n)}}g^{s,j}H_{s,i}c(\widetilde {e_i})
+\sum_{i=s=n}g^{n,j}c(\widetilde{e_n}).
\end{eqnarray*}
\noindent So for $l<n$, $\partial_{x_l}c(dx_j)(x_0)=0. $ \hfill$\Box$\\

\noindent The proof of Lemma 2.3:\\

\indent Recall, let $\nabla^L$ be the Levi-Civita connection about
$g^M$ and
$$\nabla^L_{\partial _{x_i}}\partial_{x_j}=\sum^n_{k=1}\Gamma^k_{i,j}{\partial_{x_k}},\eqno(A.1)$$
\noindent then
$$\Gamma^k_{i,j}=\frac{1}{2}g^{kl}(\partial_{x_j}g_{li}+\partial_{x_i}g_{lj}-\partial_{x_l}g_{ij}).\eqno(A.2)$$
\noindent Let
$$\partial_{x_i}=\sum_kh_{ik}\widetilde{e_k};~~\widetilde{e_i}=\sum_k\widetilde{h_{ik}}\partial_{x_k},\eqno(A.3)$$
\noindent then the matrix $[h_{ik}]$ and $[\widetilde{h_{ik}}]$
are invertible, and $\widetilde{h_{ik}}(x_0)=\delta_{ik}.$ By
(A.1) and (A.3), then
\begin{eqnarray*}
\nabla^L_{\widetilde{e_i}}\widetilde{e_t}(x_0)
&=&\nabla^L_{\partial _{x_i}}(\sum_k\widetilde{h_{tk}}\partial_{x_k})\\
&=&\sum_k\partial_{x_i}(\widetilde {h_{tk}})\partial_{x_k}
+\sum_{k,l}\widetilde{h_{t,k}}\Gamma ^l_{ik}\partial_{x_l}\\
&=&\sum_s[\sum_k\partial_{x_i}(\widetilde
{h_{tk}})h_{ks}+\sum_{k,l}\widetilde{h_{t,k}}\Gamma^l_{ik}h_{ls}]\widetilde{e_s}.
\end{eqnarray*}
\noindent By (2.7), then
$$\omega_{st}(\widetilde{e_i})(x_0)=\partial_{x_i}(\widetilde{h_{ts}})(x_0)+\Gamma_{it}^s(x_0)
=-\partial_{x_i}h_{ts}(x_0)+\Gamma_{it}^s(x_0). \eqno(A.4)$$
\noindent By (A.2) and the choices of $g^M$ and the normal coordinates of $x_0$ in $\partial M$, then we have\\

\noindent {\bf Lemma A.2}~{\it When $i<n$, then }
$$\Gamma^n_{ii}(x_0)=\frac{1}{2}h'(0);~\Gamma^i_{ni}(x_0)=-\frac{1}{2}h'(0);~\Gamma^i_{in}(x_0)=-\frac{1}{2}h'(0),$$
\noindent {\it in other cases}, $ \Gamma_{st}^i(x_0)=0.$\\
\indent By $h_{ts}=g^M(\partial_{x_t},\widetilde{e_s})=\frac{1}{\sqrt{h(x_n)}}H_{ts},~(1\leq t,s<n)$,
then we have\\

\noindent {\bf Lemma
A.3}~~$-\partial_{x_i}h_{ts}(x_0)=\frac{1}{2}h'(0)$ {\it if}
$i=n,~t=s<n$.
{\it In other cases}, $-\partial_{x_i}h_{ts}(x_0)=0.$\\
\noindent By Lemma A.2 and A.3 , (A.4), then we prove Lemma 2.3. \hfill$\Box$\\

\noindent{\bf Acknowledgement:}~~The author is indebted to Professor
Weiping Zhang for his encouragement and support. He thanks
Professors Xianzhe Dai and Siye Wu for their helpful discussions on
the gravitational action for manifolds with boundary. \\

\noindent{\bf References}\\

\noindent [A] T. Ackermann, {\it A note on the Wodzicki residue,} J.
Geom. Phys., 20, 404-406, 1996.\\
\noindent [B]  N. H. Barth, {\it The fourth-order gravitational action for manifolds with boundaries,}
Class. Quantum Grav. 2, 497-513, 1985.\\
\noindent [BGV] N. Berline, E. Getzler, M. Vergne, {\it Heat Kernals and Dirac Operators,}
Springer-Verlag, Berlin, 1992.\\
\noindent [C1] A. Connes, {\it Quantized calculus and applications,}
XIth International Congress of Mathematical Physics (paris,1994),
15-36, Internat Press, Cambridge, MA, 1995.\\
\noindent [C2] A. Connes. {\it The action functinal in
noncommutative
geometry,} Comm. Math. Phys., 117:673-683, 1998.\\
\noindent [FGLS] B. V. Fedosov, F. Golse, E. Leichtnam, and E.
Schrohe. {\it The noncommutative residue for manifolds with
boundary,} J. Funct.
Anal, 142:1-31,1996.\\
\noindent [FGV] H. Figueroa, J. Gracia-Bond\'{i}a, and J.
V\'{a}rilly, {\it Elements of Noncommutative Geometry,} Birkh\"{a}user Boston 2001.\\
\noindent [GS] G. Grubb, E. Schrohe, {\it Trace expansions and the
noncommutative residue for manifolds with boundary,} J. Reine Angew.
Math., 536:167-207, 2001.\\
 \noindent [Gu] V.W. Guillemin, {\it A new proof of Weyl's
formula on the asymptotic distribution of eigenvalues}, Adv. Math.
55 no.2, 131-160, 1985.\\
\noindent [H] S. W. Hawking, {\it General Relativity. An Einstein Centenary Survey,} Edited by S. W. Hawking and W.
Israel, Cambridge University Press,Cambridge-New York, 1979.\\
\noindent [K] D. Kastler, {\it The Dirac operator and gravitiation,}
Commun. Math. Phys, 166:633-643, 1995.\\
\noindent [KW] W. Kalau and M.Walze, {\it Gravity, non-commutative
geometry, and the Wodzicki residue,} J. Geom. Phys., 16:327-344, 1995.\\
\noindent [RS] S. Rempel and B. W. Schulze, Index theory of elliptic
boundary problems, Akademieverlag, Berlin, 1982.\\
 \noindent [S] E.
Schrohe, {\it Noncommutative residue, Dixmier's trace, and
heat trace expansions on manifolds with boundary,} Contemp. Math. 242, 161-186, 1999.\\
\noindent [U] W. J. Ugalde, {\it Differential forms and the Wodzicki
residue,} arXiv: Math, DG/0211361.\\
\noindent [Wa1] Y. Wang, {\it Differential forms and the Wodzicki
residue for manifolds with boundary,} J. Geom. Phys.,
56:731-753, 2006.\\
 \noindent
[Wa2] Y. Wang, {\it Differential forms and the noncommutative
residue for manifolds with boundary in the non-product Case,} Lett.
math. Phys., 77:41-51, 2006.\\
\noindent [Wo] M. Wodzicki,  {\it Local invariants of spectral
asymmetry}, Invent.Math. 75 no.1 143-178, 1984.\\
\noindent [Y] Y. Yu, {\it The Index Theorem and The Heat Equation Method}, Nankai Tracts in Mathematics - Vol. 2, World Scientific Publishing, 2001.\\

\end{document}